\documentclass[a4paper,11pt]{amsart}
\usepackage{amssymb}
\author[I.~V.~Arzhantsev]{Ivan V.~Arzhantsev$^*$}
\thanks{$*$ Supported by GK 02.445.11.7407 (Russia)}
\address{
Department of Higher Algebra\\ Faculty of Mechanics and
Mathematics\\ Moscow State University\\ Leninskie Gory, GSP-2, Moscow, 119992, Russia}
\email{arjantse@mccme.ru}
\urladdr{http://mech.math.msu.su/department/algebra/staff/arzhan.htm}
\author[A.~P.~Petravchuk]{Anatoliy P. Petravchuk}
\address{
Algebra Department \\ Faculty of Mechanics and Mathematics\\ Kyiv
Taras Shevchenko University\\ Volodymyrskaia street, 01033 Kyiv,
Ukraine} \email{aptr@univ.kiev.ua}
\title[Closed and Irreducible Polynomials in Several Variables]
{Closed and Irreducible Polynomials in Several Variables}
\date{May 18, 2007}

\newcommand{\kk}{\Bbbk}
\newcommand{\cchar}{\mathop{\mathrm{char}}}
\newcommand{\tr}{\mathop{\mathrm{tr}}}
\newcommand{\ZZ}{\mathbb{Z}}
\newcommand{\RR}{\mathbb{R}}
\newcommand{\QQ}{\mathbb{Q}}

\newcommand{\CC}{\mathbb{C}}
\newcommand{\Ker}{\mathop{\mathrm{Ker}}}
\newcommand{\GL}{\rm{GL}} 

\newtheorem{theorem}{Theorem}
\newtheorem{corollary}{Corollary}
\newtheorem{proposition}{Proposition}
\newtheorem{lemma}{Lemma}
\theoremstyle{definition}
\newtheorem{definition}{Definition}
\newtheorem{example}{Example}

\theoremstyle{remark}
\newtheorem{remark}{Remark}
\theoremstyle{problem}

\def\cal#1{\mathcal #1}
\def\dd#1#2{\frac{\partial #1}{\partial#2} }

\def\char{{\rm char\ }}

\let\geq\geqslant

\begin{document}

\sloppy

\begin{abstract}
New and old results on closed polynomials, i.e., such polynomials
$f\in\kk[x_1,\dots,x_n]\setminus\kk$ that the subalgebra $\kk[f]$
is integrally closed in $\kk[x_1,\dots,x_n]$, are collected in the
paper. Using some properties of closed polynomials we prove the
following factorization theorem: Let
$f\in\kk[x_1,\dots,x_n]\setminus\kk$, where $\kk$ is algebraically
closed. Then for all but finite number $\mu\in\kk$ the polynomial
$f+\mu$ can be decomposed into a product $f+\mu =\alpha\cdot
f_{1\mu}\cdot f_{2\mu}\cdots f_{k\mu}, \
\alpha\in\kk^{\times}, k\geq 1$, of irreducible polynomials
$f_{i\mu}$ of the same degree $d$ not depending on $\mu$ and such
that $f_{i\mu}-f_{j\mu}\in\kk, \quad i,j=1,\ldots,k.$ An algorithm
for finding of a generative polynomial of a given polynomial $f$,
which is a closed polynomial $h$ with $f=F(h)$ for some
$F(t)\in\kk[t]$, is given. Some types of saturated subalgebras
$A\subset\kk[x_1,\dots,x_n]$ are considered, i.e., such that for
any $f\in A\setminus\kk$ a generative polynomial of $f$ is
contained in $A$.
\end{abstract}

\maketitle

\section{Introduction}

Recall that a polynomial $f\in\kk[x_1,\dots,x_n]\setminus\kk$ is called
{\it closed} if the subalgebra $\kk[f]$ is integrally closed in 
$\kk[x_1,\dots,x_n]$. It turns out that a polynomial $f$ is closed if
and only if $f$ is {\it non-composite}, i.e., 
cannot be presented in the form $f=F(g)$ for some 
$g\in\kk[x_1,\dots,x_n]$ and $F(t)\in\kk[t], \ \deg(F)>1.$ 
Because any polynomial in $n$ variables can be obtained
from a closed polynomial by taking a polynomial in one variable
from it, the problem of studying closed polynomials is of
interest. Besides, closed polynomials in two variables appear in
a natural way as generators of rings of constants of non-zero derivations.

In this paper we present some characterizations and properties of closed polynomials.
They allow us to obtain a "generic decomposition" of a polynomial 
$f\in\kk[x_1,\dots,x_n]$ over an algebraically closed field $\kk$ (Theorem~\ref{tA}),
which may be considered as an analogue of the Fundamental
Theorem of Algebra for polynomials in many variables.

Let us go briefly through the content of the paper. 
In Section~\ref{s3} we collect numerous characterizations of closed polynomials (Theorem~\ref{ttpp}).
A major part of these characterizations is contained in the union of~\cite{Ayad}, \cite{no}, \cite{Now},
\cite{pi}, \cite{sc}, etc, but some results seem to be new. In particular, implication (i) $\Rightarrow$ (iv) in
Theorem~\ref{ttpp} over any perfect field and Proposition~\ref{prmn} solve a problem stated in~\cite[Sec.~8]{Ayad}.
For all implications in Theorem~\ref{ttpp} we give complete and elementary proofs. 

Define a generative polynomial $h$ of a polynomial $f\in\kk[x_1,\dots,x_n]\setminus\kk$ as a closed
polynomial such that $f=F(h)$ for some $F\in\kk[t]$. Clearly, a generative polynomial exists for any $f$.
Theorem~\ref{ttpp} implies that a generative polynomial is unique up to affine transformations. 

The above-mentioned results allow us to prove that over an algebraically closed field $\kk$
for any $f\in\kk[x_1,\dots,x_n]\setminus\kk$ and
for all but finite number $\mu\in\kk$ the polynomial $f+\mu$ 
can be decomposed into a product 
$f+\mu =\alpha\cdot f_{1\mu}\cdot f_{2\mu}\cdots
f_{k\mu}, \ \alpha\in\kk^{\times}, k\geq 1$,
of irreducible polynomials $f_{i\mu}$ of the same degree $d$
not depending on $\mu$ and such that
$f_{i\mu}-f_{j\mu}\in\kk, \ i,j=1,\ldots,k$.
Moreover, Stein-Lorenzini-Najib's Inequality (Theorem~\ref{sl}) implies that the number of "exceptional"
values of $\mu$ is less then $\deg(f)$. The same inequality gives an estimate of the number of irreducible
factors in $f+\mu$ for exceptional $\mu$, see Theorem~\ref{tA}.  

Section~\ref{s5} is devoted to saturated subalgebras $A\subset\kk[x_1,\dots,x_n]$, 
i.e., such that for any $f\in A\setminus\kk$ a generative polynomial of $f$ is contained in $A$. 
Clearly, any subalgebra that is integrally closed in $\kk[x_1,\dots,x_n]$ is saturated. On the other hand, 
it is known that for monomial subalgebras these two conditions are equivalent. In Theorem~\ref{prp} we
characterize subalgebras of invariants $A=\kk[x_1,\dots,x_n]^G$, where $G$ is a finite group acting
linearly on $\kk[x_1,\dots,x_n]$, with $A$ being saturated. This result provides many examples
of saturated homogeneous subalgebras that are not integrally closed in $\kk[x_1,\dots,x_n]$.

In the last section we give an algorithm for finding a generative polynomial of a given polynomial $f$. 
In particular, this provides a test to determine whether 
a given polynomial is closed.
 
%%%%%%%%%%%%%%%%%%%%%%%%%%%%%%%%%%%%%%%%%%%%%%%%%%%%%%%%%%%%%%%%%%%%%%%%%%%

\section{Characterizations of closed polynomials}\label{s3}

Let $\kk$ be a field. 

\begin{definition}
A polynomial $f(x_1,\dots,x_n)\in\kk[x_1,\dots,x_n]\setminus\kk$ is called
{\it closed} if the subalgebra $\kk[f]$ is integrally closed in $\kk[x_1,\dots,x_n]$.
\end{definition}

\begin{definition}
A polynomial $f(x_1,\dots,x_n)\in\kk[x_1,\dots,x_n]\setminus\kk$ is said to be
{\it non-composite} if the condition

\smallskip

\ (*)\ \ \ \  there exist $F(t)\in\kk[t]$ and $g(x_1,\dots,x_n)\in\kk[x_1,\dots,x_n]$
such that

\ \ \ \ \ \ \ \ $f(x_1,\dots,x_n)=F(g(x_1,\dots,x_n))$

\smallskip

\noindent implies $\deg(F)=1$.
\end{definition}

Without loss of generality, we assume below that the leading coefficient of $f$ (say, with respect to the lexicographic
order, $x_1>\dots>x_n$) equals 1. 

\smallskip

Let $\cal M$ be the set of all subalgebras $\kk[f], \ f\in\kk[x_1,\dots,x_n]\setminus\kk$, partially ordered
by inclusion.

%%%%%%%%%%%%%%%%%%%%%%%%%%%%%%%%%%%%%%%%%%%%%%%%%%

\begin{theorem}\label{ttpp}
The following conditions on a polynomial $f\in\kk[x_1,\dots,x_n]\setminus\kk$ are equivalent:

\begin{enumerate}

\item[(i)] $f$ is non-composite;

\item[(ii)] $\kk[f]$ is a maximal element of $\cal M$;

\item[(iii)] $f$ is closed;

\item[(iv)] ($\kk$ is a perfect field) \ $f+\lambda$ is irreducible over $\overline{\kk}$ for all but finitely many 
$\lambda\in\overline{\kk}$;

\item[(v)] ($\kk$ is a perfect field) \ there exists $\lambda\in\overline{\kk}$ such that $f+\lambda$ is irreducible over $\overline{\kk}$;

\item[(vi)] ($\char\kk=0$)  there exists a (finite) family of derivations $\{D_i\}$ of the algebra $\kk[x_1,\dots,x_n]$ such that 
                             $\kk[f]=\cap_i\Ker D_i$.

\end{enumerate}

\end{theorem}

%%%%%%%%%%%%%%%%%%%%%%%%%%%%%%%%%%%%%%%%%%%%%%%%%%%%%

\begin{proof}

Implications (iii) $\Rightarrow$ (ii) $\Rightarrow$ (i) \ are obvious.

\medskip

(i) $\Rightarrow$ (iii). \ It is sufficient to prove the following

\begin{proposition}\label{prv}
For any $f\in\kk[x_1,\dots,x_n]\setminus\kk$, the integral closure $A$
of $\kk[f]$ in $\kk[x_1,\dots,x_n]$ has the form $A=\kk[h]$ for some $h\in\kk[x_1,\dots,x_n]$. 
\end{proposition}

\begin{proof}
Here we follow~\cite{ea} and~\cite{za}.
Since any non-zero prime ideal of $\kk[f]$ is maximal, this also holds for $A$.

\begin{lemma}\label{llal}
The $\kk$-algebra $A$ may be realized as a $\kk$-subalgebra of a polynomial algebra in one variable.
\end{lemma}
 
\begin{proof}
By definition, $A\subset\kk[x_1,\dots,x_n]$. If $n=1$, then we are done. If $n>1$, take $r\gg 0$ such that there exists $a\in A$ with
$a+\lambda\notin I:=(x_n+x_1^r)\lhd\kk[x_1,\dots,x_n]$ for any $\lambda\in\kk$. Here $I$ is a prime ideal of $\kk[x_1,\dots,x_n]$, 
$I\cap A$ is a prime ideal of $A$, and since the image of $A$ in $\kk[x_1,\dots,x_n]/I\cong\kk[\widehat{x_1},\dots,\widehat{x_{n-1}}]$
is not a field, the ideal $I\cap A$ is zero. This shows that $A$ is embedded into $\kk[\widehat{x_1},\dots,\widehat{x_{n-1}}]$, and we
may proceed by induction on $n$.
\end{proof}

Note that since $\kk[x_1,\dots,x_n]$ is integrally closed, the subalgebra $A$ is integrally closed too. 

\begin{lemma}\label{llbl}
Let $A$ be an integrally closed subalgebra of $\kk[y]$. Then $A=\kk[u]$ for some $u\in\kk[y]$. 
\end{lemma}

\begin{proof}
Let us assume that $A\ne\kk$. Let $K$ be the quotient field of $A$. The algebra $\kk[y]$ is integral over $A$ and 
$K\cap\kk[y]=A$. Consider an element $u=u(y)\in A$ of the smallest possible positive degree $m=\deg u(y)$. Take a variable $t$ over
$K$ and consider a polynomial $F(y,t)=u(t)-u(y)\in K[t]$.  If $F(y,t)=p(y,t)q(y,t)$, then we may assume that the
highest coefficients of $p(y,t)$ and $q(y,t)$ are in $\kk$.  Moreover, $p(y,t)$ and $q(y,t)$ are in $\kk[y,t]$, because one
may consider $F=pq$ as a decomposition in $\kk(y)[t]$. This shows that the coefficients of $p$ and $q$ are in $A$. 

Fix the lexicographic monomial order with $y>t$. Then the leading term of $F(y,t)$ does not depend on $t$. So this is the
case for $p$ and $q$. By minimality of $\deg u(y)$, either $p$ or $q$ is a constant. This proves that $F(y,t)$ is irreducible
in $K[t]$, and $[\kk(y):\kk(u)]=m$. On the other hand, if $[\kk(y):K]=d$, then a linear combination of $y^d,\dots,y,1$
is contained in $K\cap\kk[y]=A$, so $d\ge m$. Hence $K=\kk(u)$ and $\kk[u]\subseteq A\subset\kk(u)$. If
$\frac{F_1(u)}{F_2(u)}\in A$ with coprime polynomials $F_1, F_2$, then there exist polynomials $F_3, F_4$ such
that $F_1F_3+F_2F_4=1$, and $F_3\frac{F_1}{F_2}+F_4=\frac{1}{F_2}\in A$. But $A$ contains no non-constant invertible elements.
This proves that $A=\kk[u]$.
\end{proof}
\end{proof}

\medskip

(iv) $\Rightarrow$ (v) \ is obvious, since $\overline{\kk}$ is infinite.

\medskip

(v) $\Rightarrow$ (i). \ If $f=F(h)$ and  
$F(t)+\lambda=(t-\alpha_1)\dots(t-\alpha_k)$ with $\alpha_1,\dots,\alpha_k\in\overline{\kk}$, 
then $f+\lambda=(h-\alpha_1)\dots(h-\alpha_k)$. 

\medskip

(i) $\Rightarrow$ (iv). \ Let us assume that $\kk=\overline{\kk}$.  Consider a morphism 
$\phi:\kk^n\to\kk^1$, $\phi(x_1,\dots,x_n)=f(x_1,\dots,x_n)$. We are going to prove that all fibers
of this morphism except for finitely many are irreducible. It follows from the next theorem
(see, for example,~\cite[p.~139]{sh}).

\medskip

{\bf The first Bertini theorem.} \ \textit{Let $X$ and $Y$ be irreducible algebraic varieties and 
$\phi:X\to Y$ be a dominant morphism. Suppose that the subfield $\phi^*(\kk(Y))$ is algebraically
closed in $\kk(X)$.  Then there exists a non-empty open subset $U\subset Y$ such that all fibers $\phi^{-1}(y)$
over $y\in U$ are irreducible}.

\medskip

So we need only to prove the following lemma.

\begin{lemma}
Let $f\in\kk[x_1,\dots,x_n]\setminus\kk$. The subfield $\kk(f)$ is algebraically closed in $\kk(x_1,\dots,x_n)$
if and only if $f$ is a closed polynomial.
\end{lemma}

\begin{proof}
If $f=F(h)$, then $h$ is in the algebraic closure of $\kk(f)$. Conversely, suppose that
for an element $g\in\kk(x_1,\dots,x_n)$ one has an equation
$$
  b_0g^m+b_1g^{m-1}+\dots+b_m=0
$$  
with $b_i\in\kk[f]$. Then $b_0g$ is integral over $\kk[f]$, and, since $f$ is closed, one has $b_0g\in\kk[f]$. This shows that
$g\in\kk(f)$.   
\end{proof}

Finally, if $\kk$ is a non-closed perfect field, the next proposition shows that $f\in\kk[x_1,\dots,x_n]$ is
closed over $\kk$ implies that $f$ is closed over $\overline{\kk}$.

\begin{proposition}\label{prmn}
Let $f\in\kk[x_1,\dots,x_n]\setminus \kk$ and $\kk\subset L$ be a separable extension of fields. Then $f$ is closed
over $\kk$ if and only if $f$ is closed over $L$.
\end{proposition}

\begin{proof}
If $f=F(h)$ over $\kk$, then the same decomposition holds over $L$.

 Now assume that $f$ is closed over $\kk$. Consider an element $g\in L[x_1,\dots,x_n]$ integral over $L[f]$. We
shall prove that $g\in L[f]$. Since the number of non-zero coefficients of $g$ is finite, we may assume
that $L$ is a finitely generated extension of $\kk$. Then there exists a finite separable transcendence basis
of $L$ over $\kk$, i.e., a finite set $\{\xi_1,\dots,\xi_m\}$ of elements in $L$ that are algebraically independent over $\kk$
and $L$ is a finite separable algebraic extension of $L_1=\kk(\xi_1,\dots\xi_m)$. 

Let us show that $f$ is closed over $L_1$. The subalgebra $\kk[f][\xi_1,\dots,\xi_m]$ is integrally closed in
$\kk[x_1,\dots,x_n][\xi_1,\dots,\xi_m]$ \cite[Ch.~V.1, Prop.~12]{bou}. Let $T$ be the set of all non-zero elements
of $\kk[\xi_1,\dots,\xi_m]$. Then the localization $T^{-1}\kk[f][\xi_1,\dots,\xi_m]$ is integrally closed in 
$T^{-1}\kk[x_1,\dots,x_n][\xi_1,\dots\xi_m]$ \cite[Ch.~V.1, Prop.~16]{bou}. This proves that $L_1[f]$ is integrally
closed in $L_1[x_1,\dots,x_n]$.

Fix a basis $\{\omega_1,\dots,\omega_k\}$ of $L$ over $L_1$. With any element $l\in L$ one may associate a $L_1$-linear
operator $M(l): L\to L$, $M(l)(\omega)=l\omega$. Let $\tr(l)$ be the trace of this operator. It is known that 
there exists a basis $\{\omega_1^{\star},\dots,\omega_k^{\star}\}$ of $L$ over $L_1$ such that 
$\tr(\omega_i\omega_j^{\star})=\delta_{ij}$ \cite[Ch.~V.1.6]{bou}. Assume that $g=\sum_i \omega_i a_i$ with 
$a_i\in L_1[x_1,\dots,x_n]$. Any $\omega_j^{\star}$ is integral over $L_1$ and thus over $L_1[f]$. This shows that
$g\omega_j^{\star}$ is integral over $L_1[f]$. Set $K=L_1(x_1,\dots,x_n)$. The element $g\omega_j^{\star}$
determines a $K$-linear map $L\otimes_K K\to L\otimes_K K$, $b\to g\omega_j^{\star}b$. Since $g\omega_j^{\star}$
is integral over $L_1[f]$, the trace of this $K$-linear operator is also integral over $L_1[f]$ \cite[Ch.~V.1.6]{bou}.  
Note that $\tr(g\omega_j^{\star})=\sum_i a_i\tr(\omega_i\omega_j^{\star})$. On the other hand, the elements
$\{\omega_1\otimes 1,\dots,\omega_k\otimes 1\}$ form a basis of $L\otimes_K K$ over $K$. 
Hence $\tr(\omega_i\omega_j^{\star})=\delta_{ij}$ and $\tr(g\omega_j^{\star})=a_j$ is integral over $L_1[f]$. 
This shows that $a_j\in L_1[f]$ for any $j$ and thus $g\in L[f]$. 
\end{proof}

\begin{example}\label{ex2}
If the field $\kk$ is not perfect, then we can not guarantee that a 
polynomial $f$ which closed over $\kk$, will be closed over
$\overline\kk$ as well. Indeed,
let $F=\kk(\eta)$ with $\eta\notin \kk$, $\eta^p\in\kk$. The polynomial
$f(x_1,x_2)=x_1^p+\eta^p x_2^p$ is closed over $\kk$. However,
one has a decomposition $f=(x_1+\eta x_2)^p$ over $F$.
The same example works for (i) $\not\Rightarrow$ (iv) in this case. 
\end{example}

\medskip

(vi) $\Rightarrow$ (iii). \ It is easy to check that for any derivation $D:\kk[x_1,\dots,x_n]\to\kk[x_1,\dots,x_n]$ 
the kernel $\Ker D$ is integrally closed in $\kk[x_1,\dots,x_n]$, and so is the intersection of kernels.  

\medskip

(i) $\Rightarrow$ (vi). \ For any $1\le i<j\le n$ consider a derivation 
$$
D_{ij}=\dd{f}{x_i}\dd{}{x_j}-\dd{f}{x_j}\dd{}{x_i}.
$$ 
Clearly, $f\in\Ker(D_{ij})$. The following Lemmas show that for a closed $f$ one has $\kk[f]=\cap_{1\le i<j\le n} D_{ij}$.

\begin{lemma}\label{lemmaB}
Assume that $\cchar\kk=0$. 
Polynomials $f,g\in\kk[x_1,\dots,x_n]\setminus\kk$ are
algebraically dependent (over $\kk$) if and only if the rank of their
Jacoby matrix $$J(f, g)= \left (
  \begin{array}{ccc}
     \dd{f}{x_{1}} & \cdots & \dd{f}{x_{n}} \\
     \dd{g}{x_{1}} & \cdots & \dd{g}{x_{n}}
     \end{array}
\right ) $$ equals to $1,$  i.e., $$\left |
  \begin{array}{cc}
     \dd{f}{x_{i}} &  \dd{f}{x_{j}} \\
     \dd{g}{x_{i}} &  \dd{g}{x_{j}}
     \end{array}
\right | =0, \quad i,j=1, \ldots , n.$$
\end{lemma}

\begin{proof}
See \cite[Ch.~3, Th.~III]{HP} or \cite[Cor.~2]{SU}.
\end{proof}

\begin{lemma}\label{lere}
Let $\kk$ be a field. 
Polynomials $f,g\in\kk[x_1,\dots,x_n]\setminus\kk$ are
algebraically dependent (over $\kk$) if and only if there exists a
closed polynomial $h\in\kk[x_1,\dots,x_n]$ such that $f,g\in\kk[h]$.
\end{lemma}

\begin{proof}
Assume that $f, g$ are algebraically dependent. 
By the Noether Normalization Lemma, there exists an element $r\in\kk[f,g]$ such that
$\kk[r]\subset\kk[f,g]$ is an integral extension. By Proposition~\ref{prv}, 
the integral closure of $\kk[r]$ in $\kk[x_1,\dots,x_n]$ 
has a form $\kk[h]$ for some closed polynomial $h$. 

Conversely, if $f, g\in k[h]$
then these polynomials are obviously algebraically dependent.
\end{proof}

This completes the proof of Theorem~\ref{ttpp}. 
\end{proof}

%%%%%%%%%%%%%%%%%%%%%%%%%%%%%%%%%%%%%%%%%%%%%%%%%%%%%

\begin{remark}
Condition (vi) may be strengthened.
It is proved in~\cite{no2} (see also~\cite{Bo}) that in characteristic zero
for any family of $\kk$-derivations $\{D_i\}$ of a finitely generated $\kk$-algebra $A$ without zero divisors 
there exists a $\kk$-derivation $D$ of $A$ with $\Ker D=\cap_i\Ker D_i$.
\end{remark}

\begin{definition}
Let $f\in\kk[x_1,\dots,x_n]\setminus\kk$. A closed polynomial $h\in\kk[x_1,\dots,x_n]$
is called {\it a generative polynomial} of $f$, if there exists $F\in\kk[t]$ such that 
$f=F(h)$.
\end{definition}

\begin{corollary}\label{ccdd}
Let $f\in\kk[x_1,\dots,x_n]\setminus\kk$. The integral closure of the subalgebra $\kk[f]$ in
$\kk[x_1,\dots,x_n]$ coincides with $\kk[h]$, where $h$ is a generative polynomial of $f$. In particular, 
a generative polynomial of $f$ exists and is unique up to affine transformations.
\end{corollary}

\begin{remark}
If one assumes that a generative polynomial $h$ of $f$ satisfies $h(0,\dots,0)=0$ and the
leading coefficient of $h$ equals 1, then such a generative polynomial exists and is unique.
Below we shall consider generative polynomials under these assumptions.
\end{remark}

%%%%%%%%%%%%%%%%%%%%%%%%%%%%%%%%%%%%%%%%%%%%%%%%%%%%%%

We refer to~\cite{Ayad} for numerous examples and interesting classes of closed polynomials.

\begin{remark}
\begin{enumerate}
\item[(1)]~Let $B$ be an affine $\kk$-algebra and again $\cal{M}(B)$ be the set of all subalgebras $\kk[b], \ b\in B\setminus\kk$, partially ordered
by inclusion. It follows from the proof above that if $B$ is integrally closed and there exists an embedding of $B$ into a polynomial
algebra, then any element $b\in B\setminus\kk$ is contained in a unique maximal element of $\cal{M}(B)$. 
For example, consider $B=\kk[x,y]/(y^2-x^3-x-1)$. Here $y^2\in\kk[x]\cap\kk[y]$, thus $B$ can not be embedded into 
a polynomial algebra.

\item[(2)] Let us recall that a subalgebra $B$ in $A$ is called {\it algebraically closed} if any
element of $A$ that is algebraic over $B$ is contained in $B$. It follows from Lemma~\ref{lere}
that a subalgebra $\kk[f]\subset\kk[x_1,\dots,x_n]$ is algebraically closed if and only
if it is integrally closed. For an arbitrary subalgebra this is no longer true (take
$\kk[x,xy]\subset\kk[x,y]$). 

\item[(3)] If a subalgebra $B\subset\kk[x_1,\dots,x_n]$ is 2-generated, then one can not
bound the number of generators of its integral closure $A$ in $\kk[x_1,\dots,x_n]$. For example, it is
easy to check that for $B=\kk[x_1,x_1x_2^m]$ the minimal number of generators of
$A=\kk[x_1,x_1x_2,\dots,x_1x_2^m]$ is $m+1$.
\end{enumerate}
\end{remark}

Finally, let us give one more observation following from Theorem~\ref{ttpp}. 
Denote by $P$ the set of all irreducible polynomials from
$\kk[x_1,\dots,x_n]$ and consider the binary
relation on this set: $f\simeq g \quad \Longleftrightarrow \exists c_{1}\in
\kk^{\times}, c_{2}\in\kk: f=c_{1}g+c_{2}.$ It is easy to see that
$\simeq$ is an equivalence relation. Choose arbitrarily a
polynomial from each equivalence class and denote by $\pi$ the set
of such polynomials.

\begin{corollary}\label{ccll}
There is a decomposition $\kk[x_1,\dots,x_n]=\cup _{p\in \pi}\kk[p],$ where
$\kk[p_i]\cap \kk[p_j]=\kk$ for all $p_{i}, p_{j}\in \pi$, $p_i\ne p_j$.
\end{corollary}

The set of nonconstant polynomials in one variable $\kk[t]\setminus
\kk$ forms a semigroup relatively to composition of polynomials
$(f\circ g)(t)=f(g(t))$. This semigroup acts naturally on
the set of nonconstant polynomials $\kk[x_1,\dots,x_n]\setminus\kk$ 
by the rule $f(h)=f(h(x_1,\dots,x_n))$
for any $h\in\kk[x_1,\dots,x_n]$. Corollary~\ref{ccll}
shows that the set $\kk[x_1,\dots,x_n]\setminus\kk$ can be
partitioned into disjoint union of "orbits" relatively to this action
and each orbit contains at least one  initial element (such that
every element of the orbit is its image but this element has no
preimages of smaller degree). Besides, this initial element is
determined up to affine transformations and can be chosen
irreducible.

%%%%%%%%%%%%%%%%%%%%%%%%%%%%%%%%%%%%%%%%%%%%%%%%%%%%%%%%%%%%%%%%%%%%%%%%

\section{A factorization theorem}\label{s4}

Let us assume in this section that our ground field $\kk$ is algebraically closed.
Theorem~\ref{ttpp} states that for a closed polynomial $h\in\kk[x_1,\dots,x_n]$ 
the polynomial $h+\lambda$ may be reducible only for finitely many $\lambda\in\kk$.
Denote by $E(h)$ the set of $\lambda\in\kk$ such that $h+\lambda$ is reducible 
and by $e(h)$ the cardinality of this set. Stein's inequality claims that
$$
e(h)<\deg f.
$$ 
Now for any $\lambda\in\kk$ consider a decomposition 
$$
h+\lambda=\prod_{i=1}^{n(\lambda,h)} h_{\lambda, i}^{d_{\lambda,i}}
$$
with
$h_{\lambda, i}$ being irreducible. A more precise version of Stein's inequality
is given in the next theorem.

\begin{theorem}[Stein-Lorenzini-Najib's Inequality]\label{sl}  
Let $h\in\kk[x_1,\dots,x_n]$ be a closed polynomial. Then
$$
\sum_{\lambda} (n(\lambda,h)-1) < \min_{\lambda}(\sum_i \deg(h_{\lambda,i})).
$$
\end{theorem}

This inequality has rather long history. Stein~\cite{st} proved his inequality 
in characteristic zero for $n=2$. For any $n$ over $\kk=\CC$ this inequality
was proved in~\cite{Cygan}. In 1993, Lorenzini~\cite{lo} obtained the inequality 
as in Theorem~\ref{sl} in any characteristic, but only for $n=2$ 
(see also \cite{ka} and \cite{vi}). Finally, in~\cite{naj} the proof for an
arbitrary $n$ was reduced to the case $n=2$.

Now take any $f\in\kk[x_1,\dots,x_n]\setminus\kk$, $\mu\in\kk$ and consider a decomposition
$$
f+\mu=\alpha\cdot\prod_{i=1}^{n(\mu,f)} f_{\mu, i}^{d_{\mu,i}}
$$
with $\alpha\in\kk^{\times}$ and $f_{\mu,i}$ being irreducible.

\medskip

Let us state the main result of this section.

\begin{theorem}\label{tA}
Let $f\in\kk[x_1,\dots,x_n]\setminus\kk$. There exists a finite subset 
$E(f)=\{\mu_1,\dots,\mu_{e(f)} \mid \mu_i\in\kk \}$ with $e(f)<\deg f$ such that

\begin{enumerate}
\item for any $\mu\notin E(f)$ one has $f+\mu=\alpha\cdot f_{\mu,1}\cdot f_{\mu,2}\cdot\cdot\cdot f_{\mu,k}$,
where all $f_{\mu,i}$ are irreducible and $f_{\mu,i}-f_{\mu,j}\in\kk$; 

\smallskip

\item $f_{\mu,i}-f_{\nu,j}\in\kk^{\times}$ for any $\mu,\nu\notin E(f)$ with $\nu\ne\mu$;
in particular, the degree $d=\deg(f_{\mu,i})$ does not depend on $i$ and $\mu$;

\smallskip

\item  $\deg(f_{\mu,i})\le d$ for any $\mu\in\kk$;

\smallskip

\item 
$$
\sum_{\mu} (n(\mu,f)-\frac{\deg(f)}{d}) < \min_{\mu} (\sum_{i=1}^{n(\mu,f)} \deg(f_{\mu,i})).
$$
\end{enumerate}
\end{theorem}

\begin{proof}
Let $h$ be the generative polynomial of $f$ and $f=F(h)$. Then 
$$
F(h)+\mu=\alpha\cdot (h+\lambda_{\mu,1})\cdot\cdot\cdot (h+\lambda_{\mu,k})
$$
for some $\lambda_{\mu,1},\dots,\lambda_{\mu,k}\in\kk$. Hence for any $\mu$ with $\lambda_{\mu,1},\dots,\lambda_{\mu,k}\notin E(h)$ 
we have a decomposition of $f+\mu$ as in (1). Note that $\lambda_{\mu,i}\ne\lambda_{\nu,j}$ for $\mu\ne\nu$. 
This proves (2) with $d=\deg(h)$ and gives the inequalities 
$$
e(f) \le e(h) < \deg(h) \le \deg(f). 
$$

Any $f_{\mu,i}$ is a divisor of some $h+\lambda$. This implies (3).

\smallskip

Finally, (4) may be obtained as:
$$
\sum_{\mu} (n(\mu,f)-\frac{\deg(f)}{d}) \le \sum_{\lambda} (n(\lambda,h)-1) < 
$$
$$
< \min_{\lambda} (\sum_i \deg(h_{\lambda,i}))
\le \min_{\mu} (\sum_j \deg(f_{\mu,j})).
$$
\end{proof}

\begin{remark}
It follows from the proof of Theorem~\ref{tA} that 
\begin{itemize}
\item $E(f)=\{ -F(-\lambda) \mid \lambda\in E(h)\}$;
\item if $f$ is not closed, then $e(f)<\frac{1}{2}\deg(f)$. 
\end{itemize}
\end{remark} 
 
\begin{example}
Take $f(x_1,x_2)=x_1^2x_2^4-2x_1^2x_2^3+x_1^2x_2^2+2x_1x_2^3-2x_1x_2^2+x_2^2+1$.

\smallskip

Here $h=x_1x_2(x_2-1)+x_2$ and $F(t)=t^2+1$. It is easy to check that $E(h)=\{0,-1\}$, thus $E(f)=\{-1,-2\}$.
We have decompositions:

\medskip

$\mu=-1: \ \ f-1=x_2^2(x_1x_2-x_1+1)^2;$

\smallskip

$\mu=-2: \ \ f-2=(x_2-1)(x_1x_2+1)(x_1x_2(x_2-1)+x_2+1);$

\smallskip

$\mu\ne -1, -2: \ \ f+\mu=(x_1x_2(x_2-1)+x_2+\lambda)(x_1x_2(x_2-1)+x_2-\lambda),$

$\ \ \ \ \ \ \ \ \ \ \ \ \ \ \ \ \ \ \  \lambda^2=-1-\mu.$

\medskip

In this case $\deg(f)=6$, $d=3$, $\sum_{\mu} (n(\mu,f)-2)=1$ and
$$
\min_{\mu} (\sum_i \deg(f_{\mu,i})) = \min\{3,6,6\}=3.
$$
\end{example} 
 
%%%%%%%%%%%%%%%%%%%%%%%%%%%%%%%%%%%%%%%%%%%%%%%%%%%%%%%%%%%%%%%%%%%%%%%

\section{Saturated subalgebras and invariants of finite groups}\label{s5}

Let $\kk$ be a field.

\begin{definition}
A subalgebra $A\subseteq\kk[x_1,\dots,x_n]$ is said to be {\it
saturated} if for any $f\in A\setminus\kk$ the generative polynomial of $f$ is
contained in $A$.
\end{definition}

Clearly, the intersection of a family of saturated subalgebras in
$\kk[x_1,\dots,x_n]$ is again a saturated subalgebra. So we may define {\it the
saturation} $S(A)$ of a subalgebra $A$ as the minimal saturated
subalgebra containing $A$.

If $A$ is integrally closed in $\kk[x_1,\dots,x_n]$, then $A$ is
saturated. By Theorem~\ref{ttpp}, if $A=\kk[f]$, then the converse is
true. Moreover, the converse is true if $A$ is a monomial
subalgebra. In order to prove it, consider a submonoid $P(A)$
in $\ZZ_{\ge 0}^n$ consisting of multidegrees of all monomials in $A$.
Then monomials corresponding to elements of the "saturated"
semigroup $P'(A)=(\QQ_{\ge 0}P(A))\cap\ZZ_{\ge 0}^n$ are generative elements
of $A$. On the other hand, it is a basic fact of toric geometry
that the monomial subalgebra corresponding to $P'(A)$ is
integrally closed in $\kk[x_1,\dots,x_n]$, see for
example~\cite[Sec.~2.1]{ful}.

Now we come from monomial to homogeneous saturated subalgebras.
The degree of monomials $\deg(\alpha x_1^{i_1}\dots x_n^{i_n})=i_1+\dots i_n$ defines a
$\ZZ_{\ge 0}$-grading on the polynomial algebra $\kk[x_1,\dots,x_n]$.
Recall that a subalgebra $A\subset\kk[x_1,\dots,x_n]$ is called {\it homogeneous} if
for any element $a\in A$ all its homogeneous components belong
to $A$.

Consider a subgroup $G\subset \GL_n(\kk)$. The linear action
$G:\kk[x_1,\dots,x_n]$ determines the
homogeneous subalgebra $\kk[x_1,\dots,x_n]^G$ of $G$-invariant
polynomials.

\begin{theorem}\label{prp}
Let $G\subseteq\GL_n(\kk)$ be a finite subgroup.
Then the subalgebra $A=\kk[x_1,\dots,x_n]^G$ is saturated in $\kk[x_1,\dots,x_n]$ if
and only if $G$ admits no non-trivial homomorphisms $G\to\kk^{\times}$.
\end{theorem}

\begin{proof}
Assume that there is a non-trivial homomorphism $\phi:G\to\kk^{\times}$.
Let $G_{\phi}$ be the kernel of $\phi$ and $G^{\phi}=G/G_{\phi}$.
Then $G^{\phi}$ is a finite cyclic group of some order $k$ and it may
be identified with a subgroup of $\kk^{\times}$.

\begin{lemma}\label{llmm}
Let $H$ be a cyclic subgroup of order $k$ in $\kk^{\times}$.
Then any finite dimensional (over $\kk$) $H$-module $W$ is a
direct sum of one-dimensional submodules.
\end{lemma}

\begin{proof}
The polynomial $X^k-1$ annihilates the linear operator $P$ in
$\GL(W)$ corresponding to a generator of $H$. By assumption, $X^k-1$ is a
product of $k$ non-proportional linear factors in $\kk[X]$.
This shows that the operator $P$ is diagonalizable.
\end{proof}

\begin{lemma}
Let $H\subset G$ be a proper subgroup. Then
$\kk[x_1,\dots,x_n]^H\ne\kk[x_1,\dots,x_n]^G$.
\end{lemma}

\begin{proof}
Let $K$ be a field and $G$ a finite group of its
automorphisms. By Artin's Theorem~\cite[Ch.VI, Th.1.8]{la},
$K^G\subset K$ is a Galois extension and  \\
$[K:K^G]=|G|$.
This implies $\kk(x_1,\dots,x_n)^H\ne\kk(x_1,\dots,x_n)^G$.
The implication
$$
 \frac{f}{h}\in\kk(x_1,\dots,x_n)^G \ \Longrightarrow \
 \frac{f\prod_{g\in G, g\ne e} g\cdot f}{h\prod_{g\in G, g\ne e} g\cdot f}\in\kk(x_1,\dots,x_n)^G
$$
shows that $\kk(x_1,\dots,x_n)^G$ (resp. $\kk(x_1,\dots,x_n)^H$) is the quotient field of
$\kk[x_1,\dots,x_n]^G$ (resp. $\kk[x_1,\dots,x_n]^H$), thus $\kk[x_1,\dots,x_n]^H\ne\kk[x_1,\dots, x_n]^G$.
\end{proof}

 Now we may take a finite-dimensional $G$-submodule $W\subset \kk[x_1,\dots,x_n]^{G_{\phi}}$
which is not contained in $\kk[x_1,\dots,x_n]^G$. Then $W$ is a $G^{\phi}$-module.
By Lemma~\ref{llmm}, one may find a $G^{\phi}$-eigenvector $h\in W$,
$h\notin\kk[x_1,\dots,x_n]^G$. Then $h^k\in\kk[x_1,\dots,x_n]^G$ and
$\kk[x_1,\dots,x_n]^G$ is not saturated.

\smallskip

Conversely, assume that any homomorphism $\chi:G\to\kk$ is
trivial. If $h$ is a generative element of a polynomial
$f\in\kk[x_1,\dots,x_n]^G$, then for any $g\in G$ the element
$g\cdot h$ is also a generative element of $f$. By
Corollary~\ref{ccdd}, the generative element is unique up to
affine transformation. Without loss of generality we can  assume
that the constant term of $h$ is zero. Then the element $g\cdot h$
has obviously zero constant term and  by Corollary~\ref{ccdd} this
element  is proportional to $h$ for any $g\in G$. Thus $G$ acts on
the line $\langle h\rangle$ via some character. But any character
of $G$ is trivial, so $h\in\kk[x_1,\dots,x_n]^G$, and
$\kk[x_1,\dots,x_n]^G$ is saturated.
\end{proof}

\begin{remark}
Since all coefficients of the polynomial
$$
F_f(T)=\prod_{g\in G} (T-g\cdot f)
$$
are in $\kk[x_1,\dots,x_n]^G$, any element $f\in\kk[x_1,\dots,x_n]$ is integral over
$\kk[x_1,\dots,x_n]^G$. Thus Theorem~\ref{prp} provides
many saturated homogeneous subalgebras that are not integrally
closed in $\kk[x_1,\dots,x_n]$.
\end{remark}

\begin{corollary}
Assume that $\kk$ is algebraically closed and $\cchar\kk=0$.
\begin{enumerate}
\item The subalgebra  $\kk[x_1,\dots,x_n]^G$ is saturated in
$\kk[x_1,\dots,x_n]$ if and only if $G$ coincides with its
commutant.
\item The saturation of $\kk[x_1,\dots,x_n]^G$ is $\kk[x_1\dots,x_n]$
if and only if $G$ is solvable.
\end{enumerate}
\end{corollary}

\begin{example}
In general, the saturation $S(A)$ is not generated by generative elements
of elements of $A$. Indeed, take any field $\kk$ that contains a primitive
root of unit of degree six. Let $G=S_3$ be the permutation group acting
naturally on $\kk[x_1,x_2,x_3]$ and $A_3\subset S_3$ be the alternating
subgroup. The proof of Theorem~\ref{prp} shows that
any generative element of an $S_3$-invariant is an $S_3$-semiinvariant and thus belongs to
$\kk[x_1,x_2,x_3]^{A_3}$. On the other hand, $S(\kk[x_1,x_2,x_3]^{S_3})=\kk[x_1,x_2,x_3]$.
\end{example}

\begin{example}
It follows from Theorem~\ref{prp} that the property of a subalgebra to be saturated is not preserved
under field extensions. Let us give an explicit example of this effect.

Let $\kk=\RR$ and $G$ be the cyclic group of order three acting on $\RR^2$ by rotations.
We begin with calculation of generators of the algebra of invariants $\RR[x,y]^G$.
Consider the complex polynomial algebra $\CC[x,y]=\RR[x,y]\oplus \text{i}\RR[x,y]$
with the natural $G$-action. Then $\CC[x,y]^G=\RR[x,y]^G\oplus \text{i}\RR[x,y]^G$.
Put $z=x+\text{i}y$, $\overline{z}=x-\text{i}y$. Clearly, $\CC[x,y]=\CC[z,\overline{z}]$,
and $G$ acts on $z, \overline{z}$ as $z\to \epsilon z$,
$\overline{z}\to\overline{\epsilon}\overline{z}$, where $\epsilon^3=1$.
This implies
$\CC[z,\overline{z}]^G=\CC[f_1,f_2,f_3]$ with $f_1=z^3$,
$f_2=\overline{z}^3$ and $f_3=z\overline{z}$. Finally,
$\RR[x,y]^G=\RR[\text{Re}(f_i), \text{Im}(f_i); \  i=1,2,3]=
\RR[x^3-3xy^2, y^3-3x^2y, x^2+y^2]$.

By Theorem~\ref{prp}, the subalgebra $\RR[x,y]^G$
is saturated in $\RR[x,y]$. On the other hand, the subalgebra $\CC[x^3-3xy^2, y^3-3x^2y, x^2+y^2]$
contains $x^3-3xy^2+\text{i}(y^3-3x^2y)=(x-\text{i}y)^3$.
\end{example}
    
%%%%%%%%%%%%%%%%%%%%%%%%%%%%%%%%%%%%%%%%%%%%%%%%%%%%%

\section{An algorithmic approach to closed polynomials}\label{s2}

Let us assume in this section that $\char\kk=0$.
Fix a homogeneous monomial order $\succ$ on the set $M(n)$ of monomials in
$x_1,\dots,x_n$, i.e., a total order satisfying the following conditions:

\begin{enumerate}
\item[(i)] $m\succ 1$ for any $m\in M(n)\setminus\{1\}$;

\item[(ii)] $m_1\succ m_2$ implies $mm_1\succ mm_2$ for all $m,m_1,m_2\in M(n)$;

\item[(iii)] if $\deg m_1>\deg m_2$, then $m_1\succ m_2$. 
\end{enumerate}

\smallskip

Take a polynomial $f\in\kk[x_1,\dots,x_n]\setminus\kk]$ such that the leading coefficient of $f$ (with
respect to $\succ$) is 1 and $f(0,\dots,0)=0$. Let $h$ be the generative polynomial of $f$ satisfying
the same assumptions. One has $f=F(h)$ with $F$ having the highest coefficient equals 1 and $F(0)=0$.

\medskip

Define {\it the multiplicity} of a monomial $m=x_1^{i_1}\dots
x_n^{i_n}$ as $d(m)=GCD(i_1,\dots,i_n)$. If $f=F(g)$,
$\deg(F)=k$, and $\overline{f}$ is the leading monomial of $f$, 
then $d(\overline{f})$ is divisible by $k$. In particular, if
$d(\overline{f})=1$, then $f$ is closed.

\medskip

Below we give an algorithm computing the generative polynomial of
a given polynomial $f(x_1,\dots,x_n)$.

\smallskip

{\bf Step 1.}\ Find a sequence $D(f)=(d_1>\dots>d_s)$ of all divisors
of $d(\overline{f})$ greater than 1. Put $j=1$.

\smallskip

{\bf Step 2.}\ Take the divisor $k=d_j$. Put $m_1=x_1^{i_1'}\dots
x_n^{i_n'}$, where $i_p'=\frac{i_p}{k}$. Let $m_1\succ
m_2\succ\dots\succ m_N\succ 1$ be the set of all monomials that do
not exceed $m_1$. Consider $h=m_1+\alpha_2m_2+\dots+\alpha_Nm_N$
with indeterminate coefficients $\alpha_j$. Let us find $\alpha_j$
inductively. Suppose that the coefficient in $f$ of the monomial
$m_1^{k-1}m_j$ equals $b_j$. Then $\alpha_2=\frac{b_2}{k}$ and,
for $j=3,\dots,N$, $\alpha_j=\frac{b_j-K_j}{k}$, where $K_j$ is
the coefficient of $m_1^{k-1}m_j$ in
$(m_1+\alpha_2m_2+\dots+\alpha_{j-1}m_{j-1})^k$.

\smallskip

{\bf Step 3.}\ Put $F(t)=t^k+\beta_1t^{k-1}+\dots+\beta_{k-1}t$
with indeterminate coefficients $\beta_l$. Assuming that $f=F(h)$
with $h$ found at Step 2, one may calculate $\beta_l$ inductively
(from $l=1$ to $k-1$) looking at the coefficient of
$m_1^{k-l}$.

\smallskip

{\bf Step 4.}\ Check the equality $f=F(h)$ with $h$ and $F$ found
at Steps 2 and 3 respectively. If the equality holds, then $h$
is the generative polynomial for $f$. (Indeed, if $h=F_1(h_1)$
with $\deg F_1>1$, then $f=F(F_1(h_1))$ and this expression
corresponds to the divisor $\deg F(F_1(t))>k$ of
$d(\overline{f})$, a contradiction.) If the equality does not hold
and $j<s$, then put $j:=j+1$ and go to Step 2. If the equality
does not hold and $j=s$, then $h=f$ is the generative polynomial
of $f$.

\smallskip

\begin{example}\label{ex1}
Consider the case $n=2$ and take the homogeneous lexicographic order with
$x_1\succ x_2$. For $f(x_1,x_2)=x_1^4+2x_1^2x_2+x_2^2$ one has
$\overline{f}=x_1^4$, $d(m)=4$ and $D(f)=(4, 2)$.

1)\ $k=d_1=4$. Here $m_1=x_1$, $m_2=x_2$ and $N=2$. One easily find that
$\alpha_2=0$, hence $h=x_1$. Moreover,
$F(t)=t^4+\beta_1t^3+\beta_2t^2+\beta_3t$. Fulfilling Step 3, we get
$\beta_1=\beta_2=\beta_3=0$. So, the equality $f=F(h)=x_1^4$ does not
hold.

2)\ $k=d_2=2$. Here $m_1=x_1^2$ and
$$
h=x_1^2+\alpha_2x_1x_2+\alpha_3x_2^2+\alpha_4x_1+\alpha_5x_2.
$$
Using Step 2, we get $\alpha_2=\alpha_3=\alpha_4=0$ and $\alpha_2=1$.
Put $F(t)=t^2+\beta_1t$, we find $\beta_1=0$ and the equality
$f=F(h)=(x_1^2+x_2)^2$ holds. So the generative polynomial for
$f$ is $h=x_1^2+x_2$.
\end{example}

\begin{remark}
If $\cchar\kk=p$ and $GCD(p,d(\overline{f}))=1$, then our algorithm also works.
\end{remark}

Finally, let us present some observations that may speed up the
above algorithm. With any term $\alpha x_1^{i_1}\dots x_n^{i_n}$,
$\alpha\ne 0$, of the polynomial $f$ one may associate an integral
point $(i_1,\dots,i_n)$ in the real vector space $\RR^n$. Let
$N(f)$ be the Newton polytope of $f$, i.e., the convex hull of all
$(i_1,\dots,i_n)$. Consider the set $V$ of all vertices of $N(f)$.
Let $$
 V_0=\{(i_1,\dots,i_n)\in V \mid \forall \ (j_1,\dots,j_n)\ne(i_1,\dots,i_n)\in V
\ \exists \ s \ : \ i_s>j_s \}.
$$

\begin{lemma}
The elements of $V_0$ parametrize terms of $f$ that are the leading terms with
respect to some monomial order on $M(n)$.
\end{lemma}

\begin{proof}
If $j_s\ge i_s$ for any $s$, then $x_1^{j_1}\dots x_n^{j_n}\succ
x_1^{i_1}\dots x_n^{i_n}$ for any monomial order $\succ$.
Conversely, if $(i_1,\dots,i_n)\in V_0$, then
$N(f)-(i_1,\dots,i_n)$ intersects the positive octant of $\RR^n$
only at zero. Hence there is a linear function
$l(z)=\omega_1z_1+\dots+\omega_nz_n$ such that the hyperplane
$l(z)=0$ separates the positive octant and $N(f)-(i_1,\dots,i_n)$,
and all $\omega_i$ are positive. (Here $z_i$ are coordinates in
$\RR^n$.) Moreover, one may assume that $\omega_i$ are linearly
independent over $\QQ$.

With any sequence $\Omega=(\omega_1,\dots,\omega_n)$ of positive real numbers
that are linearly independent over $\QQ$ one may associate a monomial
(non-homogeneous) order $\succ_{\Omega}$ defined as
$$
 x_1^{i_1}\dots x_n^{i_n} \succ_{\Omega} x_1^{j_1}\dots x_n^{j_n}
\Leftrightarrow \omega_1 i_1+\dots+\omega_n i_n>\omega_1 j_1+\dots+\omega_n j_n.
$$

Clearly, $(i_1,\dots,i_n)$ represents the leading term with respect to
$\succ_{\Omega}$ defined by the coefficients of $l(z)$.
\end{proof}

Now for any $v\in V_0$, let $m_v$ be the leading monomial of $f$
with respect to the corresponding order, and $d_1(f)=GCD(d(m_v)
\mid v\in V_0)$. Considering the leading terms of $F(h)$ with
respect to all possible monomial orders, one shows that the
sequence $D(f)$ at Step 1 may be replaced by the sequence
$D_1(f)=(\tilde d_1>\dots>\tilde d_p)$ of all divisors of $d_1(f)$
greater than 1.

 Let us return to Example~\ref{ex1}. Here $V=V_0=\{(4,0),(0,2)\}$ and
$D_1(f)=(2)$. Hence the case $k=4$ may be excluded.

%%%%%%%%%%%%%%%%%%%%%%%%%%%%%%%%%%%%%%%%%%%%%%%%%%%%%%%%%%%%%%%%%%%%%%

\end{document}